\documentclass[12pt]{amsart}

\input{epsf}
\usepackage[english]{babel}
\usepackage[dvips]{graphicx}
\usepackage{amsmath,amssymb,amsfonts,latexsym}
\usepackage{graphics,psfig}

\newtheorem{theorem}{\textbf{Theorem}}[section]

\newtheorem{proposition}[theorem]{\textbf{Proposition}}

\newtheorem{example}[theorem]{\textbf{Example}}

\def\vv{{\bf v}}
\def\ww{{\bf w}}

\title{\textsc{Independent sets in certain classes of (almost)
regular graphs}}

\begin{document}

\maketitle
\thispagestyle{empty}
\begin{center}
Alexander Burstein\\
Department of Mathematics \\ Iowa State University\\
Ames, IA 50011-2064, USA\\
\texttt{burstein@math.iastate.edu}\\[1.3ex]
Sergey Kitaev\\
Department of Mathematics \\ University of Kentucky\\
Lexington, KY 40506-0027, USA\\
\texttt{kitaev@ms.uky.edu}\\[1.3ex]
Toufik Mansour\\
Department of Mathematics \\ Haifa University\\
31905 Haifa, Israel\\
\texttt{toufik@math.haifa.ac.il}
\end{center}

\begin{abstract}
We enumerate the independent sets of several classes of regular
and almost regular graphs and compute the corresponding generating
functions. We also note the relations between these graphs and
other combinatorial objects and, in some cases, construct the
corresponding bijections.

\bigskip


\bigskip

\end{abstract}


\section{Introduction} \label{section_1}

Let $I(G)$ denote the number of independent sets of a graph $G$.
This number can be determined for some special classes of graphs
(see~\cite{Sillke} for a survey). For instance, $I(G)$ was studied
for \emph{grid graphs} (see~\cite{CalkinWilf}), \emph{multipartite
complete graphs}, and \emph{path-} and \emph{cyclic-schemes}. Many
of these numbers are given by certain combinations of Fibonacci
numbers, some others by Lucas numbers.

In this paper, we study four classes of graphs. To define these
classes, we recall that a \emph{line graph} $L(G)$ of a graph $G$
is obtained by associating a vertex with each edge $G$ and
connecting two vertices with an edge if and only if the
corresponding edges of $G$ are adjacent. Also, recall that a
\emph{cycle graph} $C_{\ell}$ is a graph on $\ell$ nodes containing a
single cycle through all nodes.

We now give the definitions of our classes.

\begin{description}
\item[Class 1] Let $G_{\ell}^{1}=C_{\ell}$, the $\ell$-cycle graph.
We obtain $G_{\ell}^{2}$ by superimposing the line graph $L(G_{\ell}^{1})$
onto the graph $G_{\ell}^{1}$, that is splitting each edge of $G_{\ell}^{1}$
with the corresponding vertex of $L(G_{\ell}^{1})$ and then adding the
edges of $L(G_{\ell}^{1})$. More generally, $G_{\ell}^{n}$ is obtained by
superimposing $L^{n-1}(G_{\ell}^1)$ onto $G_{\ell}^{n-1}$. For
example, in Figure~\ref{fig1}, we have the graphs $G_3^4$ and
$G_4^3$, respectively, if one ignores the dashed edges. Clearly,
all but $n$ of the ${\ell}^{n}$ nodes of $G_{\ell}^n$ have degree 4,
and we say that $G_{\ell}^n$ is an \emph{almost} 4-regular graph.

\begin{center}
\begin{figure}[h]
\hspace*{20pt}
\epsfxsize=200.0pt
\epsffile{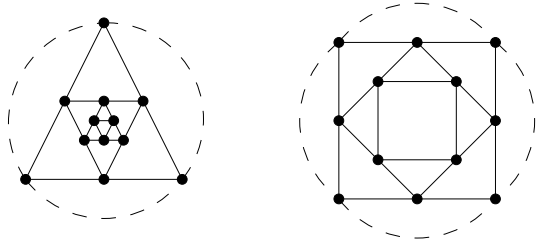}
\caption{Examples of (almost) 4-regular graphs under consideration.}
\label{fig1}
\end{figure}
\end{center}

\smallskip

\item[Class 2] The graph $R_{\ell}^n$ is obtained from $G_{\ell}^n$
by duplicating the edges of $G_{\ell}^{1}$. For example, in
Figure~\ref{fig1}, the extra edges are the dashed edges, and by
adding them we get the graphs $R_3^4$ and $R_4^3$ respectively.
So, to get $R_{\ell}^n$ we add $\ell$ additional edges to
$G_{\ell}^n$, and it is easy to see that $R_{\ell}^n$ is a
4-regular graph.

\smallskip

\item[Class 3] Let $K_{\ell}^1=K_{\ell}$, a complete graph on ${\ell}$ nodes.
Put the $\ell$ nodes of $K_{\ell}^1$ on a circle and draw the remaining
$(\ell-1)!-\ell$ edges. Call the first $\ell$ edges \emph{external} and the
remaining edges, \emph{internal}. Then construction of
$K_{\ell}^n$ is similar to that of $G_{\ell}^n$, except that:

\begin{enumerate}
\item The basis of the construction is now $K_{\ell}$, rather than
$C_{\ell}$.

\item On each iteration $i$, the graph superimposed onto
$R_{\ell}^{i-1}$ is not $L^i(K_{\ell})$ but rather the complete graph on the
nodes of $L^i(C_{\ell})$, the line graph of $C_{\ell}$ formed by the
external edges of $K_{\ell}^1$.
\end{enumerate}

In Figure~\ref{fig2}, we show how to construct $K_4^3$ from
$K_4^2$ (the dashed edges should be ignored). In that figure, the
external edges of respective complete graphs are in bold. We also
remark that the internal edges do not intersect each other.
Moreover, Figure~\ref{fig2} suggests a convenient way of
representing $K_4^3$, where each node of the graph lies only on
internal edges incident with that node. We achieve that by putting
the nodes of the line graphs under consideration off the centers
of the corresponding edges. Indeed, if we were using the centers
of the external edges, $K_4^3$ would look as in Figure~\ref{fig3},
which is misleading since, for example, the node $a$ does not lie
on the edge $bc$.

The graph $K_{\ell}^n$ is almost $(\ell+1)$-regular, since all but
$\ell$ of its ${\ell}^n$ nodes have degree $\ell+1$. Moreover, it
follows from our definitions that $K_3^n=G_3^n$.

\begin{center}
\begin{figure}[h]
\hspace*{20pt}
\epsfxsize=350.0pt
\epsffile{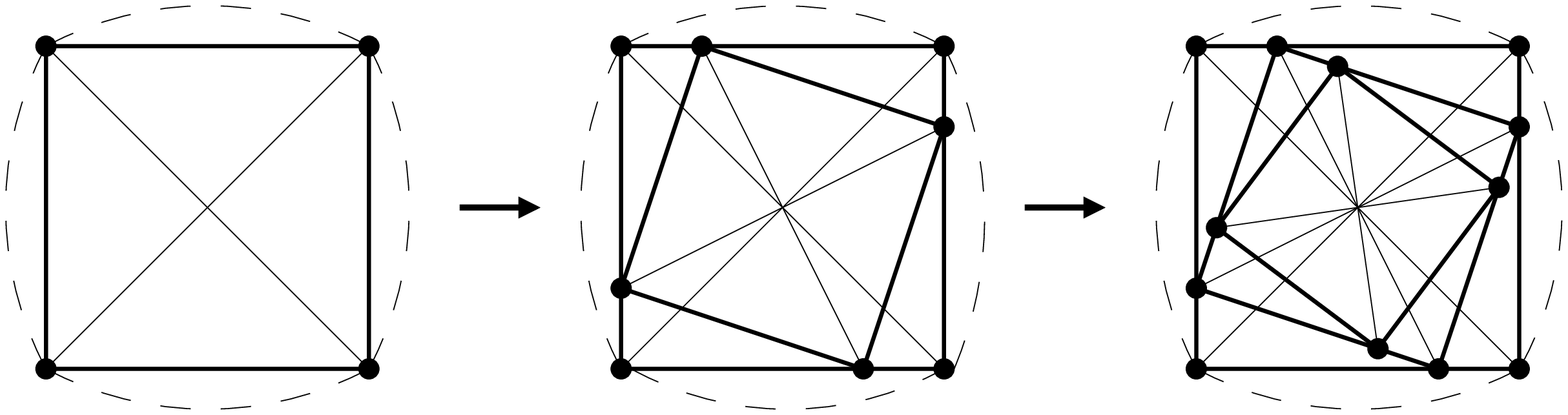}
\caption{A way of constructing $K_4^3$ and $P_4^3$.} \label{fig2}
\end{figure}
\end{center}

\smallskip

\item[Class 4]
The graph $P_{\ell}^n$ is obtained from the graph $K_{\ell}^n$ by
duplicating the external edges in the graph $K_{\ell}^{1}$. For
example, in Figure~\ref{fig2}, the extra edges are the dashed
edges, and by adding them we construct the graph $P_4^3$ from
$K_4^3$. So to get $P_{\ell}^n$ we add $\ell$ extra edges to
$K_{\ell}^n$, and it is easy to see that $P_{\ell}^n$ is an
$(\ell+1)$-regular graph.
\end{description}

\begin{center}
\begin{figure}[h]
\hspace*{20pt}
\epsfxsize=100.0pt
\epsffile{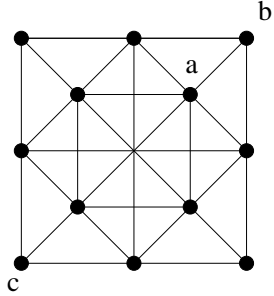}
\caption{A bad presentation of $K_4^3$.} \label{fig3}
\end{figure}
\end{center}

Let $g_{\ell}(n)$, $r_{\ell}(n)$, $k_{\ell}(n)$ and $p_{\ell}(n)$ denote
the number of independent sets in the graphs $G_{\ell}^n$,
$R_{\ell}^n$, $K_{\ell}^n$ and $P_{\ell}^n$ respectively. In our
paper we study these numbers. In Section~\ref{section_3}, we give
an algorithm for calculating all these numbers. For the numbers
$p_{\ell}(n)$, we provide an explicit generating function (see
Theorem~\ref{propp}). However, in order to illustrate our approach
to the problem, in Section~\ref{section_2} we consider $g_3(n)$
and find an explicit formula for it.

Our choice of the graphs to study was motivated by the so called
{\em de Bruijn graphs}, which are defined as follows. A de Bruijn
graph is a directed graph $\vec{G}_n=\vec{G}_n(V,E)$, where the
set of vertices $V$ is the set of all the words of length $n$ in a
finite alphabet $A$, and there is an arc from $v_i = (v_{i1},
\ldots ,v_{in})$ to $v_j = (v_{j1}, \ldots , v_{jn})$ if
\begin{center}
$v_{i2} = v_{j1}$, $v_{i3} = v_{j2}, \ldots , v_{in} =
v_{j{(n-1)}}$,
\end{center}
that is when the words $v_i$ and $v_j$ overlap by $(n-1)$ letters.

The de Bruijn graphs were first introduced (for the alphabet
$A=\{0,1\}$) by de Bruijn in 1944 for enumerating the number of
code cycles. However, these graphs proved to be a useful tool for
various problems related to the subject of combinatorics on words
(e.g. see~\cite{evdok1, evdok2, golomb}). It is known that the
graph $\vec{G}_n$ can be defined recursively as
$\vec{G}_n=L(\vec{G}_{n-1})$. The authors were interested in
studying other graphs defined recursively using the operation of
taking line graphs (with natural bases), which could give
interesting applications. Also, with our choice of graphs
($G_{\ell}^n$ and $K_{\ell}^n$), it is natural to complete them to
regular graphs ($R_{\ell}^n$ and $P_{\ell}^n$) and study these
graphs. It turns out that there are combinatorial interpretations
(relations to other combinatorial objects) for the number of
independent sets for some of our graphs, and we mention these
relations in Sections~\ref{section_2} and~\ref{section_3}.
Moreover, we construct a direct bijection describing such a
relation for $P_4^n$ (see Proposition~\ref{psets} and the
discussion that follows).

\section{The numbers $g_3(n)$.}\label{section_2}

Let us first find an explicit formula for $g_3(n)$.

It is clear that for any independent set of the graph $G_3^n$, we
can label a node of $G_3^n$ 1 if this node is in the independent
set, and label it 0 otherwise. Thus, our purpose is to count the
number of triangles having either 0 or 1 in each node and such
that no two adjacent nodes are both assigned 1s. We call such
triangles \emph{legal}.

\begin{center}
\begin{figure}[h]
\hspace*{20pt}
\epsfxsize=300.0pt
\epsffile{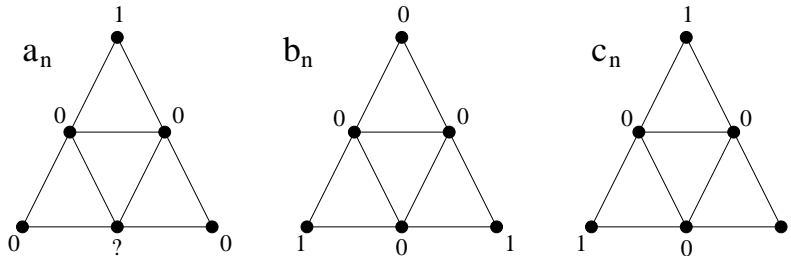}
\caption{Auxiliary parameters $a_n$, $b_n$, and $c_n$.} \label{fig4}
\end{figure}
\end{center}

In order to get a recursion for $g_3(n)$, we introduce three
auxiliary parameters $a_n$, $b_n$, and $c_n$, which are the
numbers of legal triangles that, up to rotation, have specific
numbers in the nodes of the biggest triangle (see
Figure~\ref{fig4}). Since we consider only legal graphs, the 1s in
the nodes of the biggest triangle induces 0s in certain nodes of a
smaller triangle (this 0s are shown in  Figure~\ref{fig4}).

Considering all the possibilities for the numbers of the biggest
triangle, we have that
\[
g_3(n)=g_3(n-1)+3a_n+3b_n+c_n,
\]
where $g_3(n-1)$ corresponds to all 0s, and we have multiple 3 two
times because of possible rotations. Similarly, we get that
$a_n=g_3(n-2)+a_{n-1}$, and $b_n=c_n=g_3(n-2)$. This leads to the
recursion
\begin{equation}
g_3(n)=2g_3(n-1)+6g_3(n-2)-4g_3(n-3),
\label{equa1}
\end{equation}
which, under the same initial conditions, is equivalent to the
recursion
\begin{equation}
g_3(n)=4g_3(n-1)-2g_3(n-2).
\label{equa2}
\end{equation}
We define $g_3(0)=1$, since we associate the graph $G_3^0$ with
the empty graph, in which case there is only one independent set,
the empty set. Thus,
\[
g_3(n)=\frac{1}{2\sqrt{2}}\left((2+\sqrt{2})^{n+1}-(2-\sqrt{2})^{n+1}\right),
\]
and the generating function for the numbers $g_3(n)$ is
$1/(1-4x+2x^2)$. The initial values for the numbers $g_3(n)$ are:
\[
1, 4, 14, 48, 164, 560, 1912, 6528, 22288, 76096,\ldots .
\]

Preceded by 0, the sequence $\{g_3(n)\}$ is the binomial transform
of the Pell numbers
\[
P_n=\frac{(1+\sqrt{2})^n-(1-\sqrt{2})^n}{2\sqrt{2}}
\]
(see \cite[A007070]{SloanePlouffe}). These numbers can also be
interpreted as maximum bets in a poker game (also see
\cite[A007070]{SloanePlouffe}), where the first player bets 1
dollar into a pot and the $i$th player bets the amount of the
$(i-1)$st player's bet plus the resulting amount of money in the
pot. Then the number of dollars $d_n$ in the pot after $n$ bets is
given by
\[
d_{n}=2(d_{n-1}+(d_{n-1}-d_{n-2}))=4d_{n-1}-2d_{n-2}, \quad d_0=1,
d_1=4,
\]
which yields $d_n=g_3(n)$.

We remark that it would be interesting to obtain recurrence
(\ref{equa2}) directly from the graph $G_3^{n}$, rather than via
recurrence (\ref{equa1}). Unfortunately, we were unable to do
this.

\section{An algorithm for calculating $g_{\ell}(n)$, $r_{\ell}(n)$,
$k_{\ell}(n)$ and $p_{\ell}(n)$} \label{section_3}

In this section we present an algorithm for calculating
$g_{\ell}(n)$, $r_{\ell}(n)$, $k_{\ell}(n)$ and $p_{\ell}(n)$ by
using the transfer matrix method (see
\cite[Theorem~4.7.2]{Stanley1}).

\subsection{An algorithm for calculating $g_{\ell}(n)$} \label{subsection_g}
In this section, we use the transfer matrix method to obtain an
information about the sequence of $g_{\ell}(n)$.

Similarly to Section~\ref{section_2}, for any independent set of
the graph $G_\ell^n$, consider a labeling of $G_\ell^n$, where the
nodes of of the independent set are labeled $1$ and the remaining
nodes are labeled 0. For a given graph $G_\ell^n$, we define the
\emph{$n$-th level} of $G_\ell^n$ to be $G_\ell^{n}\setminus
G_\ell^{n-1}$, which is isomorphic to $G_\ell^1$. Thus, we may
think of an independent set of the graph $G_\ell^n$ as assembled
from elements chosen on each level, making sure that when we add a
new level, we create no conflict with the previous level.

The collection $\mathcal{L}_\ell$ of possible level labelings is
the set of all $(0,1)$ $\ell$-vectors $\vv=(v_1,\dots,v_\ell)$. It
will be convenient to define $v_{\ell+1}:=v_1$. Then
$\vv=(v_1,\dots,v_\ell)$ and $\ww=(w_1,\dots,w_\ell)$ in
$\mathcal{L}_\ell$ are a possible consecutive pair of levels in an
independent set of $G_\ell^n$ (with $\ww$ following $\vv$) if and
only if
\begin{equation}\label{eqgg}
v_i=v_{i+1}=0\mbox{ or }w_i=0, \qquad \text{ where }
i=1,2,\ldots,\ell.
\end{equation}
Thus, to obtain any independent set in the graph $G_\ell^n$, we
begin with a vector of $\mathcal{L}_\ell$, then keep adjoining
each next vector $\ww\in\mathcal{L}_\ell$ so that it
satisfies~(\ref{eqgg}) together with the previously chosen vector
$\vv\in\mathcal{L}_\ell$, until $n$ vectors have been selected.

We define a matrix $G=G_\ell$, the transfer matrix of the problem,
as follows. $G$ is a $2^\ell\times 2^\ell$ matrix of $0$s and $1$s
whose rows and columns are indexed by vectors of
$\mathcal{L}_\ell$. The entry of $G$ in position $(\vv,\ww)$ is
$1$ if the ordered pair of vectors $(\vv,\ww)$
satisfies~(\ref{eqgg}), and is $0$ otherwise. $G$ depends only on
$\ell$, not on $n$. Hence, the number of independent sets of
$G_\ell^n$, $g_{\ell}(n)$, is the first entry of the vector
$G^n\cdot (u_1,\ldots,u_{2^\ell})^T$, where $u_i=1$ if there are
the $i$th vector $\vv$ in the collection $\mathcal{L}_\ell$ has no
two consecutive $1$s, even after wrapping, (i.e. $v_i+v_{i+1}\le
1$ for all $i=1,2,\ldots,\ell$), and $u_i=0$ otherwise. Hence,
\[
g_{\ell}(n)=(1,0,\ldots,0)\cdot G^n\cdot(u_1,\ldots,u_{2^\ell})^T.
\]
For instance, when $\ell=3$, the possible level vectors are
\[
(0,0,0), (0,0,1), (0,1,0), (0,1,1), (1,0,0), (1,0,1), (1,1,0),
(1,1,1),
\]
except for the last level, where we only have
\[
(0,0,0), (0,0,1), (0,1,0), (1,0,0).
\]
If we index the rows and the columns of the transfer matrix $G$ in
this order, then we get
\[
G=\begin{bmatrix}
1&1&1&1&1&1&1&1\\
1&0&0&0&1&0&0&0\\
1&1&0&0&0&0&0&0\\
1&0&0&0&0&0&0&0\\
1&0&1&0&0&0&0&0\\
1&0&0&0&0&0&0&0\\
1&0&0&0&0&0&0&0\\
1&0&0&0&0&0&0&0
\end{bmatrix}.
\]

The vector $(u_1,\ldots,u_{2^\ell})^T$ in this case is
$(1,1,1,0,1,0,0,0)$. If we now find the first entry of the vector
$(I-xG)^{-1}\cdot(1,1,1,0,1,0,0,0)^T$, where $I$ is the unit
matrix, then we get that the generating function for $g_3(n)$ is
given by $1/(1-4x+2x^2)$. We obtain the results for larger $\ell$
similarly.

\begin{theorem}
The generating functions for the numbers $g_4(n)$, $g_5(n)$ and
$g_6(n)$ are given, respectively, by
\begin{gather*}
\frac{1+4x-x^2-2x^3}{1-3x-14x^2+15x^3+7x^4},\\
\frac{(1+x)(1+5x-8x^2)}{1-5x-30x^2+69x^3+31x^4-22x^5},\\
\frac{1+10x-12x^2-50x^3+10x^4+20x^5-12x^6}{1-8x-66x^2+280x^3+178x^4-532x^5-84x^6+108x^7}.
\end{gather*}
\end{theorem}

We remark that the algorithm for finding the generating function
for $g_{\ell}(n)$ has been implemented in Maple, and yielded
explicit results for $\ell\le 6$.

\subsection{An algorithm for calculating $r_{\ell}(n)$} \label{subsection_r}
In this section we use the transfer matrix method to obtain an
information about the numbers $r_{\ell}(n)$. This case is similar
to that of $g_\ell(n)$ with some small differences.

We partition $R_\ell^n$ into levels just as in the case of
$G_\ell^n$, so the $n$-th level of $R_\ell^n$ is
$R_\ell^{n}\backslash R_\ell^{n-1}$.

The collection of possible levels $\mathcal{L}_\ell$ is the set of
all $\ell$-vectors $\vv$ of $0$s and $1$s such that there no
consecutive $1$s in $\vv$, that is, $v_i+v_{i+1}\ne 2$ (where we
define $v_{\ell+1}:=v_1$). Clearly, the set $\mathcal{L}_\ell$
contains exactly $L_{\ell}$ vectors where $L_\ell$ is the $\ell$th
Lucas number. For instance, $\mathcal{L}_3$ contains the vectors
$(0,0,0)$, $(0,0,1)$, $(0,1,0)$, and $(1,0,0)$.

The condition that vectors $\vv$ and $\ww$ in $\mathcal{L}_\ell$
are a possible consecutive pair of levels in an independent set of
$R_\ell^n$ is given by (\ref{eqgg}) just as for $G_\ell^n$. To
obtain any independent set in the graph $R_\ell^n$, we begin with
a vector of $\mathcal{L}_\ell$, then keep adjoining each next
vector $\ww\in\mathcal{L}_\ell$ so that it satisfies~(\ref{eqgg})
together with the previously chosen vector
$\vv\in\mathcal{L}_\ell$, until $n$ vectors have been selected.

We define a matrix the transfer matrix of the problem $R=R_\ell$
in the same way as $G=G_\ell$ in subsection \ref{subsection_g}.
Then $R$ is an $L_\ell\times L_\ell$ matrix, and the number of
independent sets of $R_\ell^n$, $r_{\ell}(n)$, is the first entry
of of the vector $R^n\cdot\mathbf{1}$ where
$\mathbf{1}=(1,1,\ldots,1)$. Hence,
\[
r_{\ell}(n)=(1,0,\ldots,0)\cdot R^n\cdot {\bf 1}.
\]
For instance, when $\ell=3$, the possible level vectors in an
independent set are
\[
(0,0,0), (0,0,1), (0,1,0), (1,0,0).
\]
If we indexed the rows and columns in this order, then the
transfer matrix is
\[
R=
\begin{bmatrix}
1&1&1&1\\
1&0&0&1\\
1&1&0&0\\
1&0&1&0\\
\end{bmatrix}.
\]
If we find the first entry of the vector
$(I-xG)^{-1}\cdot\mathbf{1}$, we get that the generating function
for $r_3(n)$ is given by $\frac{1+2x}{1-2x-2x^2}$. The initial
values for the numbers $r_3(n)$ are
\[
1, 4, 10, 28, 76, 208, 568, 1552, 4240, \ldots
\]
This sequence appears as A026150 in~\cite{SloanePlouffe}.


Similarly to the case $\ell=3$, we obtain the following results for
$\ell=4,5,6$.

\begin{theorem}
The generating functions for the numbers $r_4(n)$, $r_5(n)$ and
$r_6(n)$ are given, respectively, by
\begin{gather*}
\frac{1+4x-4x^2}{1-3x-4x^2+4x^3},\\
\frac{1+7x-6x^2}{1-4x-8x^2+6x^3},\\
\frac{1+12x-24x^2+8x^4}{(1-8x+4x^2+4x^3)(1+2x-2x^2)}.
\end{gather*}
\end{theorem}

We remark that the algorithm for finding the generating function
for $r_{\ell}(n)$ has been implemented in Maple, and yielded explicit
results for $\ell\leq 6$.

\subsection{An algorithm for calculating $k_{\ell}(n)$}
\label{subsection_k} In this section we use the transfer matrix
method yet again to obtain information about the sequences
$k_{\ell}(n)$. This case is also similar to that of $g_{\ell}(n)$,
so we will only sketch it briefly.

We partition $K_\ell^n$ into levels just as in the case of
$G_\ell^n$, so the $n$-th level of $K_\ell^n$ is
$K_\ell^{n}\backslash K_\ell^{n-1}$. The collection of possible
levels $\mathcal{L}_\ell$ is the set of all $\ell$-vectors $\vv$
of $0$s and $1$s. Vectors $\vv$ and $\ww$ in $\mathcal{L}_\ell$
are a possible consecutive pair of levels in an independent set of
$K_\ell^n$ if they satisfy (\ref{eqgg}).

We define the transfer matrix of the problem, $K=K_\ell$, in the
same way as $G_\ell$. Then $K$ is a $2^\ell\times 2^\ell$ matrix,
and the number of independent sets of $K_\ell^n$, $k_{\ell}(n)$,
is the first entry of of the vector $K^n\cdot
(u_1,\ldots,u_{2^\ell})^T$ where $u_i=1$ if the $i$th vector in
the collection $\mathcal{L}_\ell$ contains at most one nonzero
entry. Hence,
\[
k_{\ell}(n)=(1,0,\ldots,0)\cdot K^n\cdot
(u_1,\ldots,u_{2^\ell})^T.
\]
For instance, when $\ell=3$, the possible level vectors in an
independent set are
\[
(0,0,0), (0,0,1), (0,1,0), (0,1,1), (1,0,0), (1,0,1), (1,1,0),
(1,1,1).
\]
If we index the rows and columns in this order, then the transfer
matrix is
\[
K=
\begin{bmatrix}
1&1&1&1&1&1&1&1\\
1&0&0&0&1&0&0&0\\
1&1&0&0&0&0&0&0\\
1&0&0&0&0&0&0&0\\
1&0&1&0&0&0&0&0\\
1&0&0&0&0&0&0&0\\
1&0&0&0&0&0&0&0\\
1&0&0&0&0&0&0&0
\end{bmatrix}.
\]

If we find the first entry of the vector
$(I-xK)^{-1}\cdot(u_1,\ldots,u_{2^\ell})^T$, we get that the
generating function for $k_3(n)$ is given by
$\frac{1}{1-4x+2x^2}$. In particular, we get that $k_3(n)=g_3(n)$
which is can also be seen directly from the definitions.
Similarly, we have the following result.

\begin{theorem}
The generating functions for the numbers $k_4(n)$, $k_5(n)$ and
$k_6(n)$ are given, respectively, by
\begin{gather*}
\frac{1+2x+3x^2}{1-3x-14x^2+15x^3+7x^4},\\
\frac{1+x+12x^2-8x^3}{1-5x-30x^2+69x^3+31x^4-22x^5},\\
\frac{1-x+38x^2-72x^3-8x^4+30x^5}{1-8x-66x^2+280x^3+178x^4-532x^5-84x^6+108x^7}.
\end{gather*}
\end{theorem}

We remark that the algorithm for finding the generating function
for $k_{\ell}(n)$ has been implemented in Maple, and yielded explicit
results for $\ell\leq 6$.

\subsection{An algorithm for calculating $p_{\ell}(n)$}
In this section we use the transfer matrix method to obtain
information about the sequences $p_{\ell}(n)$.

In this section we use the transfer matrix method once more to
obtain information about the sequences $k_{\ell}(n)$. This case is
also similar to that of $g_{\ell}(n)$, so we will only sketch it
briefly.

We partition $P_\ell^n$ into levels the same way as $G_\ell^n$, so
the $n$-th level of $P_\ell^n$ is $P_\ell^{n}\backslash
P_\ell^{n-1}$. The collection of possible levels
$\mathcal{L}_\ell$ is the set of all $\ell$-vectors $\vv$ of $0$s
and $1$s. Vectors $\vv$ and $\ww$ in $\mathcal{L}_\ell$ are a
possible consecutive pair of levels in an independent set of
$P_\ell^n$ if they satisfy (\ref{eqgg}).

The collection of possible levels $\mathcal{L}_\ell$ is the set of
all $\ell$-vectors $\vv=(v_1,\dots,v_\ell)$ of $0$s and $1$s such
that $v_1+\cdots+v_\ell\le 1$. Clearly, the set $\mathcal{L}_\ell$
contains exactly $\ell+1$ vectors which are $(0,\cdots,0)$ and
$(0,\ldots,0,1,0,\ldots,0)$. The condition that vectors $\vv$ and
$\ww$ in $\mathcal{L}_\ell$ are a possible consecutive pair of
levels in an independent set of $P_\ell^n$ is given by
(\ref{eqgg}).

We define the transfer matrix of the problem, $P=P_\ell$, in the
same way as $G_\ell$, $R_\ell$ and $K_\ell$. We define a matrix
$P=P_\ell$, the transfer matrix of the problem, as follows. $P$ is
an $(\ell+1)\times (\ell+1)$ matrix of $0$s and $1$s whose rows
and columns indexed by vectors of $\mathcal{L}_\ell$. Therefore,
it is easy to see that
\[
P=
\begin{bmatrix}
1&1&1&1&1&\cdots&1&1&1\\
1&0&0&1&1&\cdots&1&1&1\\
1&1&0&0&1&\cdots&1&1&1\\
\vdots&&&&&\vdots&&&\vdots\\
1&1&1&1&1&\cdots&1&0&0\\
1&0&1&1&1&\cdots&1&1&0
\end{bmatrix}_{(\ell+1)\times(\ell+1)}.
\]
\begin{theorem}\label{propp}
The generating function for $p_{\ell}(n)$ is given by
\[
\sum_{n\geq0}p_{\ell}(n)x^n=\frac{1+2x}{1-(\ell-1)x-2x^2}.
\]
\end{theorem}
\begin{proof}
We want to find the first entry of the vector $(I-xP)^{-1}\cdot
{\bf 1}$ which means we must find the first row, say
$(e_1,\ldots,e_{\ell+1})$, of the matrix $(I-xP)^{-1}$. By solving
the system of equations
\[
(I-xP)^{-1}\cdot(e_1,\dots,e_{\ell+1})^T=(1,0,\dots,0)^T,
\]
we get that
\[
e_1=\frac{1-(\ell-2)x}{1-(\ell-1)x-2x^2}\mbox{ and
}e_j=\frac{x}{1-(\ell-1)x-2x^2}\mbox{ for }j\ge 2.
\]
Hence, the first entry of the vector $(I-xP)^{-1}\cdot{\bf 1}$ is
given by
\[
\sum_{i=1}^{\ell+1} e_i=\frac{1+2x}{1-(\ell-1)x-2x^2}.
\]
\end{proof}

For instance, when $\ell=3$, the generating function for $p_3(n)$
is given by $\frac{1+2x}{1-2x-2x^2}$. One can see, in particular,
that $p_3(n)=r_3(n)$, which follows directly from the definitions.

In the case $\ell=4$ the initial values of the numbers $p_4(n)$
are
\[
1, 5, 17, 61, 217, 773, 2753, 9805, 34921, 124373, \dots
\]

The same sequence turns out to appear in \cite{guymoser}
(see~\cite[A007483]{SloanePlouffe}). Thus, the following proposition is
true:

\begin{proposition}\label{psets}
The number of independent sets in the graph $P_4^n$ is equal to
the number of (possibly empty) subsequences of the sequence
$\{1,2,\ldots,2n+1\}$ in which each odd member has an even
neighbor.
\end{proposition}

Here, the neighbors of an integer $m$ are $m-1$ and $m+1$. In the
case $n=1$, the sequences appearing in the proposition are
$\epsilon, 2, 23, 12, 123$, where $\epsilon$ is the empty
sequence. In this case, we can find a direct bijection between the
objects in Proposition~\ref{psets} as described below.

\begin{center}
\begin{figure}[h]
\hspace*{20pt} \epsfxsize=350.0pt \epsffile{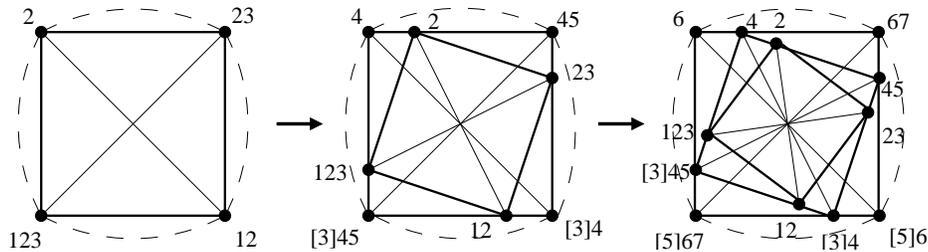}
\caption{A labeling of $P_4^3$ by subsequences of
Proposition~\ref{psets}.} \label{fig5}
\end{figure}
\end{center}

We start by labeling the vertices of the (innermost) level $n$ of
$P_4^n$ clockwise by 2, 23, 12, 123. The level $n-1$  is labeled
starting from the vertex immediately to the left of vertex labeled
2 as follows: $4,45,[3]4,[3]45$ (the meaning of brackets will be
discussed below). More generally, for $i<n$, given a level $n-i+1$
labeled clockwise with $2i,2i(2i+1),[2i-1]2i,[2i-1]2i(2i+1)$, we
label the (next outer) level $n-i$ clockwise from the inside out
with $2i+2,(2i+2)(2i+3),[2i+1](2i+2),[2i+1](2i+2)(2i+3)$ starting
from the vertex immediately to the left of vertex labeled $2i$.
(See Figure \ref{fig5} for the case $n=3$.) Each independent set
has at most one vertex on each level. Now, given any independent
set in $P_4^n$, we can write the labels of its nodes in increasing
order and delete any integer $[2i-1]$ in brackets if the sequence
also contains $2i-2$ or $2i-1$ without brackets. Erasing all
brackets now, if any, we obtain a sequence from
Proposition~\ref{psets}.

\begin{example} \label{ex:psets}
The independent sets $\{45,[5]67\}, \{4,[5]67\}, \{12,[5]6\},
\{[3]4,67\}$ correspond to the sequences $4567, 467, 1256, 3467$,
respectively.
\end{example}

For convenience, we will write down the set of nonadjacent labels
at level $n-i$ for each label at level $n-i+1$.

{\small
\begin{equation} \label{eq:rules}
\begin{split}
\text{(no label) } \epsilon&\mapsto (2i+2), (2i+2)(2i+3),
[2i+1](2i+2), [2i+1](2i+2)(2i+3)\\
2i &\mapsto [2i+1](2i+2), [2i+1](2i+2)(2i+3)\\
2i(2i+1) &\mapsto (2i+2), [2i+1](2i+2)(2i+3)\\
[2i-1]2i &\mapsto (2i+2), (2i+2)(2i+3)\\
[2i-1]2i(2i+1) &\mapsto (2i+2)(2i+3), [2i+1](2i+2)
\end{split}
\end{equation}
}

Now it is not difficult to construct an independent set given a
sequence of Proposition~\ref{psets}. We partition the sequence of
integers from 1 to $2n+1$ as follows:
\[123\,|\,45\,|\,67\,|\dots|\,2n(2n+1),\] then choose the vertices
of the independent set in the order of increasing labels using the
rules (\ref{eq:rules}). Notice that the label of the vertex at
level $n-i+1$ must contain $2i$.

\begin{example} \label{ex:pseqs}
\[
\begin{split}
4567 &\mapsto 123|\underline{45}|\underline{67} \mapsto
(\epsilon,45,[5]67) \mapsto \{45,[5]67\}\\
467 &\mapsto 123|\underline{4}5|\underline{67} \mapsto
(\epsilon,4,[5]67) \mapsto \{45,[5]67\}\\
1256 &\mapsto \underline{12}3|4\underline{5}|\underline{6}7
\mapsto (12,\epsilon,[5]6) \mapsto \{12,[5]6\}\\
\end{split}
\]
\end{example}

It can be shown that the two maps described above are inverses of
each other based on the recursive structure of sequences under
consideration.


\end{document}